\documentclass[a4paper, 12pt, 3p, times]{P-article}
\usepackage{amssymb}
\usepackage{epsfig}
\usepackage{amsfonts}
\usepackage{amsmath, array}
\usepackage{euscript}
\usepackage{amscd}
\usepackage{amsthm}
\usepackage{ulem}
\usepackage{xcolor}
\usepackage{yfonts, mathrsfs, eufrak}
\usepackage{yfonts}
\usepackage{booktabs}

\usepackage[english]{babel}
\usepackage{blindtext}
\usepackage{mathtools}
\usepackage{xcolor}
\usepackage{soul} 

\usepackage[pdftex,
pdfauthor=Prosenjit Das,
pdfsubject={Affine Algebraic Geometry, Commutative Algebra},
pdfproducer=TeXStudio,
pdfcreator=pdflatex]{hyperref}

\usepackage{hyperref}

\hypersetup{
	colorlinks=true,
}

\DeclareMathAlphabet{\mathpzc}{OT1}{pzc}{m}{it}
\newtheorem{thm}{Theorem}[section]
\newtheorem{lem}[thm]{Lemma}
\newtheorem{prop}[thm]{Proposition}

\newtheorem{cor}[thm]{Corollary}
\newtheorem{qtn}[thm]{Question}

\newdefinition{defn}[thm]{Definition}
\newdefinition{ex}[thm]{Example}
\newdefinition{rem}[thm]{Remark}
\newdefinition{note}{Note}

\newcommand{\comment}[1]{}

\newcommand\m {\mathfrak{m}}
\newcommand\p {\mathfrak{p}}
\newcommand\q {\mathfrak{q}}

\newcommand{\A}[1]{\mathbb{A}^{#1}}

\newcommand{\LND}[2]{{\rm{LND}}_{#1}(#2)}

\newcommand{\Ker}[1]{{\rm Ker}(#1)}

\newcommand{\Sym}[2]{{\rm Sym}_{#1}(#2)}

\newcommand{\ResRk}[1]{{\rm Res-Rk}(#1)}

\newcommand{\ResVarRk}[1]{{\rm ResVar-Rk}(#1)}

\newcommand{\Ht}[1]{{\rm{ht}}(#1)}

\newcommand{\Spec}[1]{\text{Spec}(#1)}
\newcommand{\bQ}{\mathbb{Q}}

\newcommand{\bR}{\mathbb{R}}
\newcommand{\bC}{\mathbb{C}}

\newcommand{\trdeg}[2]{{\rm{tr.deg}}_{#2}(#1)}

\newcommand{\grade}[1]{{\rm{grade}}\big{(}#1\big{)}}

\begin{document}
	\begin{frontmatter}
		\title{Locally nilpotent derivations on $\mathbb{A}^2$-fibrations with\\ $\mathbb{A}^1$-fibration kernels
		}
		
\author{Janaki Raman Babu}
\address{Department of Mathematics, Indian Institute of Space Science and Technology, \\
Valiamala P.O., Trivandrum 695 547, India\\
email: \texttt{janakiramanb.16@res.iist.ac.in, raman.janaki93@gmail.com}}
		
\author{Prosenjit Das\footnote{Corresponding author.}}
\address{Department of Mathematics, Indian Institute of Space Science and Technology, \\
Valiamala P.O., Trivandrum 695 547, India\\
email: \texttt{prosenjit.das@iist.ac.in, prosenjit.das@gmail.com}}
\author{Animesh Lahiri }
\address{Chennai Mathematical Institute\\
H1, SIPCOT IT Park, Siruseri\\ Kelambakkam 603 103, India\\
email: \texttt{animeshl@cmi.ac.in, 255alahiri@gmail.com}}

\begin{abstract}
In this paper, we give a characterization of locally nilpotent derivations on $\mathbb{A}^2$-fibrations over Noetherian domains containing $\mathbb{Q}$ having kernels isomorphic to $\mathbb{A}^1$-fibrations.

\smallskip
\noindent
{\small {{\bf Keywords.} Affine fibration; Locally nilpotent derivation; Kernel; Grade; Depth}

}

\end{abstract}
		
	\end{frontmatter}
	
\section{Introduction}\label{Sec_intro} \label{Sec_Intro}

\noindent
Throughout this article rings will be commutative with unity. Let $R$ be a ring and $B$ an $R$-algebra. We write $B = R^{[n]}$ to mean that $B$ is isomorphic, as an $R$-algebra, to a polynomial ring in $n$ variables over $R$. Suppose $R$ is a $\bQ$-algebra. An $R$-linear map $D:B\longrightarrow B$ is said to be an {\it $R$-derivation} if it satisfies the Leibnitz rule: $D(ab)=aD(b)+bD(a),~\forall a,b \in B$. Let $D: B \longrightarrow B$ be an $R$-derivation. The {\it kernel} of $D$, denoted by $\Ker{D}$, is defined to be the set $\{b \in B:D(b)=0\}$. $D$ is said to be a {\it locally nilpotent $R$-derivation} (abbrev. {$R$-lnd}) if, for each $b\in{B}$, there exists $n(b)\geqslant 0$ such that $D^{n(b)}(b)=0$. The set of all locally nilpotent $R$-derivations on $B$ will be denoted by $\LND{R}{B}$. Unless otherwise stated, capital letters like $X_1,\dots,X_n,Y_1,\dots,Y_m,X,Y,Z,W$ will be used to denote variables of polynomial algebras.
	
\medskip
One of the important advancements in the theory of locally nilpotent derivations, due to the works of Rentschler (\cite{Rentchler_Operations-du-Groupe}), Daigle-Freudenburg (\cite{Daigle-Freudenburg_UFD-LND-Rank-2}), Bhatwadekar-Dutta (\cite{Bhatwadekar-Dutta_LND}), Berson-van den Essen-Maubach (\cite{Essen-Maubach-Berson_Der-div-zero}) and van den Essen (\cite{Essen_Around-Cancellation}), establish that over a ring $R$ containing $\mathbb{Q}$, any fixed point free $R$-lnd $D$ (see Definition \ref{Defn_fpf}) on $R^{[2]}$ has a slice and its kernel is $R^{[1]}$.

\begin{thm} \label{Essen_FixdPtLND-R^[2]-trivKer}
Let $R$ be a ring containing $\mathbb{Q}$, $B = R[X,Y]$ and $D \in \LND{R}{B}$. If $D$ is fixed point free, then $\Ker{D}=R[f]~(=R^{[1]})$ for some $f$ in $B$, and $D$ has a slice, i.e., $B=R[f]^{[1]}$.
\end{thm} 

\smallskip\noindent
Among the above-mentioned contributions, when $R$ is a Noetherian domain, the work of Bhatwadekar-Dutta, in \cite{Bhatwadekar-Dutta_LND} draws a special attention as they have provided a necessary and sufficient condition for $\Ker{D}$ to be $R^{[1]}$. They also have described the structure of $\Ker{D}$ for {\it any} $R$-lnd $D$ under the hypothesis ``$R$ is a Noetherian normal domain".
\begin{thm} \label{BD_LND-Characterisation-Result}
Let $R$ be a Noetherian domain containing $\mathbb{Q}$, $B = R[X,Y]$ and $D \in \LND{R}{B}$. Then,
\begin{enumerate}
\item[\rm(a)] The following conditions are equivalent. 
	
\begin{enumerate}
\item[\rm(I)] $D$ is irreducible and $\Ker{D} = R^{[1]}$.
\item [\rm (II)] Either $D$ is fixed point free or $D(X)$ and $D(Y)$ form a $B$-regular sequence.
\end{enumerate}

\smallskip\noindent
If $D$ is fixed point free, then $D$ has a slice, i.e., $B = \Ker{D}^{[1]}$.

\smallskip\noindent
\item[\rm(b)] If $R$ is a normal domain which is not a field, then there exists a height one unmixed ideal $I$ of $R$ such that ${\rm Ker}(D)$ is isomorphic to the symbolic Rees-algebra ${\underset{n\geqslant 0}\bigoplus}I^{(n)}T^n$. 
\end{enumerate}
\end{thm}
\noindent
If $R$ is Noetherian, then for a non-zero $R$-lnd $D$ on $B~(= R[X,Y])$, one can easily see that $\grade{D(B)B, B}$ (see Definition \ref{grade}) is either $1$ or $2$ or $\infty$, and further one can observe that the condition (II) in Theorem \ref{BD_LND-Characterisation-Result}(a) can be replaced by $\grade{D(B)B,B}\in\{2,\ \infty\}$ (see Lemma \ref{Lem_2genIdeal_seq-grade}). Therefore, given an irreducible $R$-lnd $D$ on $B$, it follows from Theorem \ref{BD_LND-Characterisation-Result}(a) that $\Ker{D}= R^{[1]}$ if and only if $\grade{D(B)B, B} = 2 \ \text{or} \ \infty$ (see Theorem \ref{BD_LND-Characterisation-Result_Grade}).

\medskip
Recently, in \cite[Theorem 4.4]{Babu-Das_Struct_A2-fib_FPF-LND},  Babu-Das proved the following result (also see \cite[Corollary 3.4]{Kahoui_Fixd-pt-fr} and \cite[Corollary 3.2]{Kahoui_A2-fib_triviality-criterion}) which extends Theorem \ref{Essen_FixdPtLND-R^[2]-trivKer} (and hence partially Theorem \ref{BD_LND-Characterisation-Result}) to the case $B$ is an $\A{2}$-fibration over $R$ (see Definition \ref{Defn_affine fibrations}).

\begin{thm} \label{Babu-Das_A2-fib-FPF-LND}
	Let $R$ be a Noetherian ring containing $\mathbb{Q}$, $B$ an $\A{2}$-fibration over $R$ and $D \in \LND{R}{B}$. If $D$ is fixed point free, then $\Ker{D}$ is an $\A{1}$-fibration over $R$ and  $D$ has a slice, i.e., $B = \Ker{D}^{[1]}$.
\end{thm}

\noindent
In view of Theorem $\ref{Babu-Das_A2-fib-FPF-LND}$, one naturally asks whether Theorem \ref{BD_LND-Characterisation-Result} has an analogue when ``$B=R[X,Y]$'' is replaced by  ``$B$ is an $\A{2}$-fibration over $R$''.

\begin{qtn} \label{Qtn_main}
Let $R$ be a Noetherian domain containing $\mathbb{Q}$ and $B$ an $\A{2}$-fibration over $R$. 
\begin{enumerate}
\item[\rm(a)] In the spirit of Theorem \ref{BD_LND-Characterisation-Result}$({\rm a})$, is it possible to characterize the locally nilpotent $R$ derivations on $B$ having kernels isomorphic to $\A{1}$-fibrations over $R$?
\item[\rm(b)] If $R$ is a normal domain, then, in the spirit of Theorem \ref{BD_LND-Characterisation-Result}$({\rm b})$, is it possible to describe the explicit structure of $\Ker{D}$ for any $R$-lnd $D$? 
	
\end{enumerate}
\end{qtn}

\noindent
This article gives complete answer to Question \ref{Qtn_main}. 

\medskip \noindent
In Section \ref{Sec_grade}, we observe some basic facts which we use in the subsequent sections of this paper; some of them are of interest on their own. Namely, Proposition \ref{Prop_A2Fib-LND-grades} establishes that for an $\A{2}$-fibration $B$ over $R$ and a non-zero $R$-lnd $D$, grade of the ideal $D(B)B$ can not be other than $1,2,\infty$. 

\smallskip\noindent
Section \ref{Sec_main} of this article discusses our main results, and primarily focuses on part (a) of Question \ref{Qtn_main}. We establish a necessary and sufficient condition for the kernel of an $R$-lnd $D$ to be an $\A{1}$-fibration over $R$ (see Theorem \ref{Thm_A2-Fib_A1-Ker}).

\smallskip\noindent
\textbf{Theorem A1.} Let $R$ be a Noetherian domain containing $\bQ$, $B$ an $\A{2}$-fibration over $R$ and $D \in \LND{R}{B}$ be such that for each $\p\in\Spec{R}$, the induced $R_{\p}$-lnd $D_\p$ is irreducible. Then, $\Ker{D}$ is an $\A{1}$-fibration over $R$ if and only if $\grade{D(B)B, B} \in \{2, \ \infty \}$. 

\medskip\noindent To prove Theorem \ref{Thm_A2-Fib_A1-Ker}, we have proved an auxiliary result which provides a necessary and sufficient condition for the kernel of an $R$-lnd $D$ on $\Sym{R}{M}$ where $M$ is a projective $R$-module of rank two, to be $\Sym{R}{N}$ for some rank one projective $R$-module $N$ (see Proposition \ref{Prop_Triv-A2-Fib_Sym-Ker}).

\smallskip\noindent
\textbf{Theorem A2.} Let $R$ be a Noetherian domain containing $\bQ$, $B=\Sym{R}{M}$ for some rank two projective $R$-module $M$ and $D \in \LND{R}{B}$ be such that for each $\p\in\Spec{R}$, the induced $R_{\p}$-lnd $D_\p$ is irreducible. Then, $\Ker{D}$ is  isomorphic to $\Sym{R}{N}$ for some projective rank one $R$-module $N$ if and only if $\grade{D(B)B, B} \in \{2, \ \infty \}$. 

\medskip\noindent
In Section \ref{Sec_main}, we also observe some interesting corollaries to Theorem \ref{Thm_A2-Fib_A1-Ker}. Namely, Corollary \ref{Cor_A2-forms_Sym-Ker} provides a necessary and sufficient condition for the kernel of an $R$-lnd $D$ on an $\A{2}$-form over $R$ (see Definition \ref{Defn_A2forms}) to be $\Sym{R}{N}$ for some projective rank one $R$-module $N$; Corollary \ref{Cor_resrank} provides a sufficient condition for the kernel of an $R$-lnd $D$ on an $\A{n}$-fibration over $R$ to be $\A{n-1}$-fibration over $R$ when $\grade{D(B)B,B}$ is $2$ or $\infty$.

\medskip\noindent Section \ref{Sec_Normal} of this article discusses part (b) of Question \ref{Qtn_main}. If $R$ is a Noetherian normal domain and $D$ is an $R$-lnd on an $\A{2}$-fibration $B$ over $R$, then we have shown that the structure of $\Ker{D}$ is exactly same as the case of $B=R^{[2]}$. The proof goes along the line of the proof of Proposition 3.3 of \citep{Bhatwadekar-Dutta_LND}, and heavily uses the fact that any affine fibration over a Noetherian ring is a retract of a polynomial algebra (see Theorem \ref{Thm_kernel_Noeth_normal}).

\smallskip\noindent
{\bf Theorem B.} Let $R$ be a Noetherian normal domain containing $\bQ$, $B$ an $\A{2}$-fibration over $R$ and $D\in\LND{R}{B}\setminus\{0\}$. Then, $\Ker{D}$ has the structure of a graded ring $\bigoplus_{i\ge 0} A_i$ with $A_0 = R$ and for each $i \geqslant 1$, $A_i$ is a finite reflexive $R$-module of rank one. In fact, when $R$ is not a field, then there exists an ideal $I$ of unmixed height one in $R$ such that $A$ is isomorphic to the symbolic Rees algebra $\bigoplus_{n\ge 0} I^{(n)}T^{n}$. 

\medskip\noindent
 In Section \ref{Sec_Examples}, we discuss some examples. 

\section{Preliminaries}\label{Sec_Preli} \label{Sec_Prelim}
\noindent
In this section we set up notations, recall definitions and quote some results.

\medskip\noindent
{\bf{Notation:}}
Given a ring $R$ and an $R$-algebra $B$ we fix the following notation. 

\medskip
$
\begin{array}{lll}
R^*& : & \text{Group of units of $R$}.\\


\Spec{R} & : & \text{Set of all prime ideals of $R$}.\\

\Sym{R}{M} & : & \text{Symmetric algebra of an $R$-module $M$}.\\

\Omega_R(B) & : & \text{Universal module of $R$-differentials of $B$}.\\
		

B_b & : & \text{Localization of $B$ with respect to $\{1,b,b^2,\dots\}, ~~b\in B$}.\\

B_\p & : & \text{Localization of $B$ with respect to $R\setminus \p,~~\p\in\Spec{R}$}.\\

k(\p) & : & \dfrac{R_\p}{\p R_\p}, \text{ residue field of $R$ at $\p$}.\\

\end{array}
$

\begin{center}
{\bf Definitions.}
\end{center}

\begin{enumerate}
\item Let $R$ be a Noetherian ring. A sequence of elements $r_1,\dots,r_n\in R$ is said to form an $R$-regular sequence if, for each $i$, $1\leqslant i\leqslant n$, $r_i$ is a non-zerodivisor of $\dfrac{R}{(r_1,\dots,r_{i-1})R}$ and $(r_1,\dots,r_n)R\neq R$.
 \item\label{grade} Let $R$ be a Noetherian ring and $I$ an ideal of $R$. The {\it grade} of $I$ is denoted by $\grade{I, R}$ and is defined by 
$$\grade{I, R} = \left\{
\begin{array}{ll}

\text{length of maximal $R$-regular sequence contained in $I$}, & \text{if} \ I \ne R\\
\infty, & \text{if} \ I = R.
\end{array}
\right\}
$$

\item A reduced ring $R$ is called \textit{seminormal} if whenever $a^2 = b^3$ for some $a,b \in R$, then there exists $t \in R$ such that $t^3 = a$ and $t^2 =b$.
\item Let $B$ be an integral domain. A subring $A$ of $B$ is said to be \textit {inert} in $B$ if, for $a,b\in B\setminus\{0\}$, $ab\in A$ implies $a,b\in A$.

\smallskip\noindent
$\bullet$ It is easy to verify that the inertness property is preserved under localization; and an inert subring of a UFD is a UFD. 

\medskip\indent For the rest of the definitions we assume that $R$ is a ring and $B$ is an $R$-algebra.

\item $R$ is said to be a {\it retract} of $B$ if there exists an algebra homomorphism $\phi:B\longrightarrow R$ such that $\phi|_{R}=Id_R$.

\smallskip\noindent
$\bullet$ It is easy to see that if $R$ is a retract of $B$ and $S$ a multiplicatively closed subset of $R$, then $S^{-1}R$ is a retract of $S^{-1}B$.

\smallskip\noindent
$\bullet$ It is well-known that retraction of a UFD is a UFD (see p. 9 fn. in \citep{Eakin-Heinzer_A-Cncl-prob}).

\item \label{Defn_affine fibrations} $B$ is said to be an \textit{$\A{n}$-fibration over $R$}, if $B$ is finitely generated as an $R$-algebra, flat as an $R$-module and $B \otimes_R k(\p)=k(\p)^{[n]}$ for each $\p \in \Spec{R}$.
\item\label{Defn_A2forms} 
Suppose $k$ is a field of characteristic $p~(\geqslant 0)$ with algebraic closure 
$\bar{k}$ and $R$ a $k$-algebra. $B$ is said to be an {\it $\A{n}$-form over $R$} (with respect to $k$) 
if $B\otimes_{k}\bar{k}={(R\otimes_{k}{\bar{k}})}^{[n]}$.
\item  A derivation $D$ on $B$ is said to be \textit{reducible}, if there exists $b \in B\backslash B^*$ such that $D(B) \subseteq bB$. Otherwise, $D$ is called {\it irreducible}.

\item\label{Defn_fpf}
 A derivation $D$ on $B$ is said to be \textit{fixed point free} if $D(B)B=B$. 

\item  A derivation $D$ on $B$ is said to have a \textit{slice} $s \in B$, if $D(s) =1$. If $\mathbb{Q} \hookrightarrow B$ and $D$ is locally nilpotent with a slice $s$, then it is well known (\textit{slice theorem}) that $B=\Ker{D}[s]= \Ker{D}^{[1]}$ (see \cite[Corollary 1.26, p. 28]{Freudenburg-BookNew}).

\end{enumerate}	

\begin{center}
{\bf Preliminary results.}
\end{center}
\noindent
We now quote a few results for later use. The first one is by Hamann (\cite[Theorem 2.8]{Haman_Invariance}).
\begin{thm} \label{Hamann}
Let $R$ be a Noetherian ring containing $\mathbb{Q}$ and $B$ an $R$-algebra such that $B^{[m]} = R^{[m+1]}$ for some $m \in \mathbb{N}$. Then, $B = R^{[1]}$.
\end{thm}

A L\"uroth-type result by Abhyankar-Eakin-Heinzer (\cite[Proposition 4.1]{AEH_Coff}) states the following.

\begin{thm} \label{AEH_Luroth}
Let $R$ be a {\rm UFD} and $A$ be an inert $R$-subalgebra of $B=R^{[n]}$ such that ${\rm tr.deg}_{R}(A)=1$. Then, $A=R^{[1]}$.
\end{thm}

A classification of locally polynomial algebras by Bass-Connell-Wright (\cite[Theorem 4.4]{BCR_Local-Poly}) states

\begin{thm}\label{BCW_Sym}
Let $R$ be a ring and $B$ a finitely presented $R$-algebra. If $B_\m$ is a polynomial algebra over $R_\m$ for each maximal ideal $\m$ of $R$, then $B$ is $R$-isomorphic to the symmetric algebra $\Sym{R}{M}$ for some projective $R$-module $M$.
\end{thm}

The next result is by Swan (\cite[Theorem 6.1]{Swan_On-Sem}).

\begin{thm} \label{Swan_SemiNormality}
Let $R$ be a seminormal ring. Then, ${\rm Pic}(R) = {\rm Pic}(R^{[n]})$ for all $n \in \mathbb{N}$.
\end{thm}

Now, we state a result of Giral (\citep[Proposition 2.1]{Giral_subalgebra-of-polynomial-algebra}).
\begin{prop}\label{Giral}
Let $R$ be a Noetherian domain and $B$ an $R$-subalgebra of a finitely generated $R$-algebra. Then there exists $b\in B\setminus\{0\}$ such that $B_b$ is a finitely generated $R$-algebra.
\end{prop}
A result on finite generation of an algebra by Onoda (\cite[Lemma 2.14 \& Theorem 2.20]{O_Subring}) is as follows.

\begin{thm} \label{Onoda_FG}
Let $R$ be a Noetherian domain and $B$ an overdomain of $R$ such that $B_b$  is finitely generated over $R$ for some $b \in B \backslash \{0\}$. Then the following statements hold.
\begin{enumerate}
\item [\rm (I)] If $S$ is a multiplicatively closed subset of $R$ such that $B \otimes_{R} S^{-1}R$ is finitely generated over $S^{-1}R$, then there exists $s \in S$ such that $B_s$ is finitely generated over $R$.
\item [\rm (II)] $B$ is finitely generated over $R$ if and only if $B_\m$ is finitely generated over $R_\m$ for all maximal ideal $m$ of $R$.
\end{enumerate}
	
\end{thm}
Asanuma established the following structure theorem (\cite[Theorem 3.4]{Asanuma_fibre_ring}) of affine fibrations over Noetherian rings.
\begin{thm} \label{Asanuma_struct-fib-th}
Let $R$ be a Noetherian ring and $B$ an $\A{r}$-fibration over $R$. Then, $\Omega_R(B)$ is a projective $B$-module of rank $r$ and $B$ is an $R$-subalgebra (up to an isomorphism) of a polynomial ring $R^{[m]}$ for some $m \in \mathbb{N}$ such that $B^{[m]}=\Sym{R^{[m]}} {\Omega_R(B) \otimes_B R^{[m]}}$. Therefore, $B$ is a retract of $R^{[n]}$ for some $n$.
\end{thm}

\begin{cor} \label{Cor_A1-fib_Triv_DiffMod_Extended}
Let $R$ be a Noetherian ring containing $\bQ$ and $B$ an $\A{1}$-fibration over $R$. If $\Omega_R(B)$ is extended from $R$ (for example, $R$ is seminormal), then $B= {\rm Sym}_{R}(N)$ for some finitely generated rank one projective $R$-module $N$.
\end{cor}
\begin{proof}
Follows from Theorem \ref{Asanuma_struct-fib-th}, Theorem \ref{Swan_SemiNormality}, Theorem \ref{Hamann} and Theorem \ref{BCW_Sym}.
\end{proof}

\indent
Next, we quote a patching lemma by Bhatwadekar-Dutta (\citep[Lemma 3.1]{Bhatwadekar-Dutta_LND}).
\begin{lem} \label{Lem_BD-A1-patch}
Let $R$ be a Noetherian domain and $A$ an overdomain of $R$ such that $JA \cap R = J$ for every ideal $J$ of $R$. Suppose that there exist non-zero elements $x, y \in R$ satisfying the conditions:
\begin{enumerate}
	\item [\rm (i)] $x$ and $y$ form an $R$-regular sequence,
	\item [\rm (ii)] $A_x = R_x^{[1]}$ and $A_y = R_y^{[1]}$,
	\item [\rm (iii)] $A=A_x \cap A_y$.
\end{enumerate}
Then, $A$ has a graded ring structure ${\bigoplus_{i\geqslant 0}A_i}$, where $A_0=R$ and for each $i\geqslant 1$, $A_i$ is a reflexive $R$-module of rank one. In fact, $A$ is $R$-isomorphic  as a graded $R$-algebra to the symbolic Rees-algebra $\bigoplus_{n\geqslant 0}I^{(n)}T^n$ of a reflexive ideal $I$ in $R$ of height one.
\end{lem}
 
\medskip
\noindent
We now record a theorem on separable $\A{1}$-forms over rings due to Dutta (\cite[Theorem 7]{Dutta_A1-form}). 

\begin{thm}\label{AKD1}
Let $k$ be a field, $L$ a separable field extension of $k$, $R$ a $k$-algebra and 
$B$ an $R$-algebra such that $B\otimes_{k}L$ is isomorphic to the symmetric algebra 
of a finitely generated rank one projective module over $R\otimes_{k}L$. Then $B$ is isomorphic to the symmetric algebra of a finitely generated rank one projective module over $R$. 
\end{thm}

\section{Grade of the ideal $D(B)B$}\label{Sec_grade}
	
\noindent	
In this section,  we observe that the condition (II) of Theorem \ref{BD_LND-Characterisation-Result}(a) can be replaced by the equivalent condition ``$\grade{D(B)B,B}\in\{2,\infty\}$". We also prove that for an $R$-lnd $D$ on an $\A{2}$-fibration $B$ over $R$, $\grade{D(B)B,B}\in\{1,2,\infty\}$. First, we observe an easy lemma.		 
\begin{lem} \label{Lem_2genIdeal_seq-grade}
Let $R$ be a Noetherian domain and $I= (x_1, x_2)$ an ideal of $R$. Then the following are equivalent.
\begin{enumerate}
\item [\rm (I)] $x_1, x_2$ form an $R$ regular sequence.
\item [\rm (II)] ${\rm grade}(I,R) =2$.
\end{enumerate}
\end{lem}
	
\begin{proof}
{\rm(I)} $\implies$ {\rm(II)}: Suppose $x_1, x_2$ form an $R$-regular sequence. Then $\grade{I,R}\geqslant 2$. Since each element of $I$ is trivially a zero-divisor of $R/(x_1,x_2)$, we get $x_1,x_2$ is a maximal $R$-sequence in $I$. Hence $\grade{I,R}=2$.
		
\smallskip
\noindent		
{\rm(II)} $\implies$ {\rm(I)}:
Suppose, $\grade{I,R} =2$. Then $I\subsetneqq R$. Since $R$ is a domain, $x_1$ is an $R$-regular element. If possible suppose $x_2$ is a zero divisor of $R/x_1R$. Then, there exists $r\in R$ such that $rx_2\in x_1R$. Hence, for any $r_1,r_2\in R$ we have $r(r_1x_1+r_2x_2)\in x_1R$, i.e., any element of $I$ is a zero-divisor of $R/x_1R$. So, $x_1$ is a maximal $R$-regular sequence in $I$, contradicting the hypothesis $\grade{I,R}=2$. 
	\end{proof}
	
\noindent	
As an immediate consequence of Lemma \ref{Lem_2genIdeal_seq-grade} we get the following.
	
\begin{cor} \label{Cor_R2_gr2_DxDy-seq}
Let $R$ be a Noetherian domain, $B=R[X,Y]$ and $D\in{\rm LND}_{R}(B)$. Then, $D(X)$ and $D(Y)$ form a $B$-regular sequence if and only if $\grade{D(B),B} =2$.
\end{cor}

\noindent
In view of Lemma \ref{Lem_2genIdeal_seq-grade} and Corollary \ref{Cor_R2_gr2_DxDy-seq} we clearly see that Theorem \ref{BD_LND-Characterisation-Result} can be restated as follows.
	
\begin{thm} \label{BD_LND-Characterisation-Result_Grade}
Let $R$ be a Noetherian domain containing $\mathbb{Q}$, $B=R[X,Y]$ and $D\in{\rm LND}_{R}(B)$. Then the following are equivalent.
\begin{enumerate}
\item [\rm (I)] $D$ is irreducible and ${\rm Ker}(D)=R[f]~(=R^{[1]})$ for some $f \in B$.
\item [\rm (II)] $\grade{D(B),B}$ is either $2$ or $\infty$.
\end{enumerate}
\end{thm}

Now, we state and prove two easy lemmas.
\begin{lem}\label{Lem_equality of grade under ff}
  Let $R$ be a Noetherian domain and $I$ an ideal of $R$. If $B$ is faithfully flat over $R$, then $\grade{I,R}=\grade{IB,B}$.
 \end{lem}
\begin{proof}
 Since $B$ is faithfully flat over $R$, $I=IB\cap R$. Hence $I=R$ if and only if $IB=B$. So, $\grade{I,R}=\infty$ if and only if $\grade{IB,B}=\infty$.
 
 \smallskip
 \noindent
 Let $a_1,\dots,a_n$ be an $R$-regular sequence in $I$. Since $B$ is faithfully flat over $R$, $a_1,\dots,a_n$ form a $B$-regular sequence in $IB$. Suppose, $a_1,\dots,a_n$ form a maximal $R$-regular sequence in $I$. Then every element of $I$ is a zero-divisor of $\dfrac{R}{(a_1,\dots,a_n)R}$. So, for each $x\in I$ there exists $r_x\in R\setminus\{0\}$ such that $r_xx\in(a_1,\dots,a_n)R$. Let $f\in IB$. Then, there exist $x_1,\dots,x_m\in I$ and $b_1,\dots,b_m\in B$ such that $f= \displaystyle \sum_{i=1}^mx_ib_i$.
 Let $r=\displaystyle \prod_{i=1}^mr_{x_i}$. Clearly, $r\neq 0$ and $rf\in(a_1,\dots,a_n)B$. So, $f$ is a zero-divisor of $\dfrac{B}{(a_1,\dots,a_n)B}$. Thus, every element of $IB$ is a zero-divisor of $\dfrac{B}{(a_1,\dots,a_n)B}$, and hence $a_1,\dots,a_n$ is a maximal $B$-regular sequence in $IB$.
\end{proof}

\begin{rem}\label{Rem_grade under flat extension}
 From the proof of Lemma \ref{Lem_equality of grade under ff} we observe that for an ideal $I$ of a Noetherian domain $R$ and a flat ring extension $B$ over $R$, if ${\rm grade}(I,R)=n$, then ${\rm grade}(IB,B)$ is either $n$ or $\infty$.
\end{rem}

\begin{lem} \label{Lem_grade-achieved-locally-over-base-ring}
Let $R$ be a Noetherian domain containing $\bQ$, $B$ an $R$-algebra and $D \in {\rm LND}_{R}(B)$. If ${\rm grade}(D(B),B)=i~(\neq\infty)$, then there exists $\p\in{\rm Spec}(R)$ such that ${\rm grade}(D_\p(B_\p),B_\p)=i$, where $D_\p$ is the induced $R_\p$-derivation on $B_\p$.
\end{lem}
\begin{proof}
Since $\grade{D(B)B,B}=i$, there exists $P \in \Spec{B}$ such that $D(B)B \subseteq P$ and ${\rm depth}(B_{P})=i$. Let $\p=P \cap R$. Then $D(B)B\subseteq D(B_{\p})B_{\p} = D(B)B_{\p}\subseteq D(B)B_{P}\subseteq PB_P$ and hence $\grade{D(B)B,B} \leqslant \grade{D(B)B_{\p},B_{\p}}\leqslant \grade{D(B)B_P,B_P}\leqslant{\rm depth}(B_P)$. Thus, we have\\
$\grade{D(B)B, B} = \grade{D(B_{\p})B_{\p},B_{\p}} = i$.
\end{proof}
We are now ready to show that the grade of the ideal generated by the image of a non-zero $R$-lnd on an $\A{2}$-fibration over $R$ can not be other than $1, \ 2$ and $\infty$.
\begin{prop} \label{Prop_A2Fib-LND-grades}
Let $R$ be a Noetherian domain containing $\mathbb{Q}$, $B$ an $\A{2}$-fibration over $R$ and $D \in {\rm LND}_{R}(B) \setminus \{ 0 \}$. Then, ${\rm grade}(D(B)B, B) \in \{ 1, 2, \infty\}$.
\end{prop}
\begin{proof}
If $\grade{D(B)B,B}=\infty$ then we are done. So, we suppose that $\grade{D(B)B,B}=i\neq\infty$. By Lemma \ref{Lem_grade-achieved-locally-over-base-ring}, there exists $\p \in \Spec{R}$ such that $\grade{D(B)B, B} = \grade{D(B_{\p})B_{\p},B_{\p}} = i$. 
We will show that $i \in \{ 1, \ 2, \ \infty \}$.

\medskip
\noindent
First we prove the result for the case $B = \Sym{R}{M}$ for some rank two projective $R$-module $M$. In that case $B_{\p} = R_{\p}^{[2]}$. Since $D(B_{\p})B_{\p}~(\ne 0)$ is generated by two elements, we have $\grade{D(B_{\p})B_{\p},B_{\p}} \in \{1, \ 2, \ \infty \}$, i.e., $i \in \{1, \ 2, \ \infty \}$.

\smallskip\noindent
We now prove the general case. Since $B$ is an $\A{2}$-fibration over the Noetherian domain $R$, by Theorem \ref{Asanuma_struct-fib-th}  $\Omega_R(B)$ is a projective $B$-module of rank $2$ and there exists $C = R^{[n]}$ for some $n \in \mathbb{N}$ such that $B$ is an $R$-subalgebra of $C$ and $\widetilde{B} : = B \otimes_R C = \Sym{C}{\Omega_R(B) \otimes_B C}$. Let $\widetilde{D}:= D \otimes_{R} 1_C$ be the trivial extension of $D$ to $\widetilde{B}$. Clearly, $\widetilde{D}$ is a $C$-lnd of $\widetilde{B}$ and $\widetilde{D}(\widetilde{B}) \widetilde{B} = D(B)\widetilde{B}$. Since $\widetilde{B}$ is faithfully flat over $B$, by Lemma \ref{Lem_equality of grade under ff}, we have $\grade{\widetilde{D}(\widetilde{B}) \widetilde{B}} = \grade{D(B)B} =i$. Since $\Omega_R(B) \otimes_B C$ is a projective $C$-module of rank $2$, by our previous case we have $\grade{\widetilde{D}(\widetilde{B}) \widetilde{B}} \in \{ 1, \ 2, \ \infty \}$, i.e., $i \in \{ 1, \ 2, \ \infty \}$.
\end{proof}

\section{Main Results}\label{Sec_main}
\noindent

\noindent
In this section, we answer part (a) of Question \ref{Qtn_main}.	We begin with an easy lemma.
\begin{lem} \label{Lem_Poly-base-extn-Sym_implies_A1-fib}
Let $R$ be a ring, $C = R[X_1,\dots,X_n]$ and $A$ an $R$-algebra. If for some $r\geqslant 0$, $A \otimes_R C$ is an $\A{r}$-fibration over $C$, then $A$ is an $\A{r}$-fibration over $R$.
\end{lem}
\begin{proof}
Let $I=(X_1,\dots,X_n)C$. Since $A \otimes_R C$ is an $\A{r}$-fibration over $C$, $(A \otimes_R C)\otimes_C C/I$ is an $\A{r}$-fibration over $C/I$. Now, the result follows from the fact that $C/I\cong R$.	
\end{proof}		
\begin{prop} \label{Prop_Triv-A2-Fib_Sym-Ker}
Let $R$ be a Noetherian domain containing $\bQ$ and $B={\rm Sym}_{R}{(M)}$ for some finitely generated rank two projective $R$-module $M$. Suppose, $D\in{\rm LND}_{R}(B)$ and $A={\rm Ker}(D)$. Then the following are equivalent.
		
\begin{enumerate}
\item [\rm (I)] ${\rm grade}\big{(}D(B)B,B\big{)}\in\{2,\ \infty\}$.
			
\item [\rm (II)] $A={\rm Sym}_{R}{(I)}$ for some invertible ideal $I$ of $R$ and $D_{\p}$ is irreducible for each $\p\in {\rm Spec}(R)$. 
\end{enumerate}
\noindent
Moreover, when ${\rm grade}\big{(}D(B)B,B\big{)}=\infty$, we have $B = A^{[1]}$.
\end{prop}
	
\begin{proof}
{\rm(I)}$\implies${\rm(II)}: Suppose, (I) holds. Let $\p \in \Spec{R}$. Then, $B_\p=R_{\p} [X,Y]~(= {R_{\p}}^{[2]})$ for some $X,Y \in B$. Let $D_{\p}$ denote the extension of $D$ to $B_{\p}$. Then, $\Ker{D_{\p}}=A_{\p}$ and $D_{\p}(B_{p})B_{p}=D(B)B_{\p}$. Since $B_{\p}$ is a flat $B$ module, by Remark \ref{Rem_grade under flat extension}, $\grade{D_{\p}(B_{\p})B_{\p},B_{\p}}$ is either $2$ or $\infty$. By Theorem \ref{BD_LND-Characterisation-Result_Grade}, $D_{\p}$ is irreducible and $A_{\p}= R_\p^{[1]}$. Let $s\in B\setminus A$ be such that $D^2(s)=0$. Then $B_{D(s)}=A_{D(s)}[s]$. This shows that $A_{D(s)}[s]$ is a finitely generated algebra over $R$ and hence $A_{D(s)}$ is finitely generated over $R$. Now, by Theorem \ref{Onoda_FG}, we conclude that $A$ is finitely generated over $R$ and hence by Theorem \ref{BCW_Sym}, we have $A=\Sym{R}{I}$ for some invertible ideal $I$ of $R$.
		
\medskip
\noindent
{\rm(II)} $\implies$ {\rm(I)}: We assume that (II) holds. By Theorem \ref{BD_LND-Characterisation-Result_Grade}, for each $\q\in\Spec{R}$, $\grade{D(B_{\q})B_{\q},B_{\q}}$ is either $2$ or $\infty$ . If possible suppose $\grade{D(B)B,B}=i$, where $i\neq 2,\infty$.  By Lemma \ref{Lem_grade-achieved-locally-over-base-ring}, there exists $\p \in \Spec{R}$ such that $\grade{D(B_{\p})B_{\p},B_{\p}} = i$, which is a contradiction. This proves that ${\rm grade}\big{(}D(B)B,B\big{)}$ is either $2$ or $\infty$.

\medskip
\noindent
For the next part, let $\grade{D(B)B,B}=\infty$. Then $D(B)B=B$, i.e., $D$ is fixed point free. Now, by 
\cite[Proposition 4.3]{Babu-Das_Struct_A2-fib_FPF-LND} (also see \cite[Corollary 3.2]{Kahoui_A2-fib_triviality-criterion}), it follows that $B = \Ker{D}^{[1]}$.
\end{proof}
\begin{rem}\label{rem_triviality_sym}
In Proposition \ref{Prop_Triv-A2-Fib_Sym-Ker}, when $\grade{D(B),B}=\infty$, it has been proved in \citep[Proposition 4.3]{Babu-Das_Struct_A2-fib_FPF-LND} (also see \cite[Corollary 3.2]{Kahoui_A2-fib_triviality-criterion}) that $B=A^{[1]}$. 
\end{rem}	
\begin{thm} \label{Thm_A2-Fib_A1-Ker}
Let $R$ be a Noetherian domain containing $\bQ$ and $B$ an $\A{2}$-fibration over $R$. Suppose, $D\in{\rm LND}_{R}(B)$ and $A={\rm Ker}(D)$. Then the following are equivalent.
		
\begin{enumerate}
\item [\rm (I)] ${\rm grade}(D(B),B)\in\{2,\ \infty\}$.
\item [\rm (II)] $A$ is an $\A{1}$-fibration over $R$ and $D_\p$ is irreducible for each $\p\in {\rm Spec}(R)$.  
\end{enumerate}  
\end{thm}
	
\begin{proof}
(I) $\implies$ (II): Suppose, (I) holds. Since $B$ is an $\A{2}$-fibration over $R$, by Theorem \ref{Asanuma_struct-fib-th}, there exists $C:=R[X_1,\dots,X_n]$  such that $B\subseteq C$ and $\widetilde{B}:=B\otimes_R C=\Sym{C}{\Omega_R(B) \otimes_B C}$. Let $\widetilde{D}:=D\otimes_{R} {1_C}$ be the trivial extension of $D$ 
to $\widetilde{B}$, 
defined by $\widetilde{D}(b\otimes_R c)=D(b)\otimes_R c$ for all $b\in B$ and $c\in C$. Then $\Ker{\widetilde{D}}=A\otimes_RC$ and $\widetilde{D}(\widetilde{B})\widetilde{B}=D(B)\widetilde{B}$. Since $\widetilde{B}$ is faithfully flat over $B$, by Lemma \ref{Lem_equality of grade under ff}, we have $\grade{\widetilde{D}(\widetilde{B})\widetilde{B},\widetilde{B}}=\grade{D(B)B,B}$. Since $\Omega_R(B)$ is a projective $B$-module of rank $2$, it follows that $\Omega_R(B) \otimes_B C$ is a projective $C$-module of rank $2$. Now, by Proposition  \ref{Prop_Triv-A2-Fib_Sym-Ker}, we get the following.
\begin{enumerate}
 \item [\rm(i)] $A\otimes_RC= \Sym{C}{N}$ for some finitely generated rank one projective $C$-module $N$,
 \item[\rm(ii)] $\widetilde{D}_Q$ is irreducible for each $Q \in \Spec{C}$. 
\end{enumerate}

\noindent
By Lemma \ref{Lem_Poly-base-extn-Sym_implies_A1-fib}, we see that (i) implies $A$ is an $\A{1}$-fibration over $R$.

\smallskip\noindent
Let $\p\in\Spec{R}$. Since $C$ is faithfully flat over $R$, there exists $P \in \Spec{C}$ such that $\p = P \cap R$. Then $B_{\p}\subseteq \widetilde{B}_{P}$. If $D_{\p}$ is the extension of $D$ to $B_{\p}$ and $\widetilde{D}_{P}$ is the extension of $\widetilde{D}$ to $\widetilde{B}_P$, then $\widetilde{D}_P|_{B_\p}=D_{\p}$ and $D_{\p}(B_{\p})\widetilde{B}_{P}=D(B)\widetilde{B}_{P}=\widetilde{D}_{P}(\widetilde{B}_P)\widetilde{B}_{P}$. If possible suppose $D_{\p}$ is reducible. Then $D_{\p}(B_{\p})\subseteq bB_{\p}$ for some non-unit $b$ of $B_{\p}$. Hence $\widetilde{D}_{P}(\widetilde{B}_P)\subseteq b\widetilde{B}_P$. Since $B_{\p}\cap(\widetilde{B}_P)^{*}=(B_{\p})^{*}$, $b$ is a non-unit of $\widetilde{B}_P$, contradicting (ii). So, we have $D_{\p}$ is irreducible for each $\p\in\Spec{R}$. 
		
\medskip
\noindent
(II)$\implies$(I): Suppose, (II) holds. By Theorem \ref{Asanuma_struct-fib-th}, it follows that there exist $C_1 = R[X_1,\dots,X_{n_1}]$ containing $A$ and $C_2 = R[Y_1,\dots,Y_{n_2}]$ containing $B$ such that $A\otimes_{R} C_1 = \Sym{C_1}{\Omega_R(A)\otimes_A C_1}$ and $B \otimes_R C_2 = \Sym{C_2}{\Omega_R(B)\otimes_B C_2}$. Let $C := C_1 \otimes_R C_2$. Then $A\otimes_{R} C = \Sym{C}{\Omega_R(A)\otimes_A C}$ and $B\otimes_R C = \Sym{C}{\Omega_R(B)\otimes_B C}$. Let $\widetilde{D}$ be the trivial extension of $D$ to $B \otimes_R C$. Let $P\in\Spec{C}$ and $\p=P\cap R$. Then $C_P$ is faithfully flat over $R_{\p}$ and hence $\widetilde{B}_P~(=C_P\otimes_{R_{\p}}B_{\p})$ is faithfully flat over $B_{\p}$. If possible suppose $\widetilde{D}_P$ is reducible. Then $\widetilde{D}(\widetilde{B})\widetilde{B}_P\subseteq b\widetilde{B}_P$, where $b\in\widetilde{B}$ is a non-unit of $\widetilde{B}_P$. Since $D_{\p}(B_{\p})\widetilde{B}_{P}~(=D(B)\widetilde{B}_{P})=\widetilde{D}_{P}(\widetilde{B}_P)\widetilde{B}_{P}$, we have $D_{\p}(B_{\p})\widetilde{B}_P\subseteq b\widetilde{B}_P$ and hence $D_{\p}(B_{\p})B_{\p}\subseteq b\widetilde{B}_P\cap{B_{\p}}$. Since $D_{\p}\neq0$, $b\widetilde{B}_P\cap{B_{\p}}\neq (0)$. Then there exists $f\in\widetilde{B}$, $g\in C\setminus P$, $a\in B$ and $r\in R\setminus \p$ such that $rbf=ag$. 

\smallskip\noindent
{\bf Case 1:} Suppose, $rf\in C$. Then $rf\in C\setminus P$ and hence $b~(=a\dfrac{g}{rf})\in aB_P$. Since $\dfrac{g}{rf}\in(\widetilde B_P)^*$, $b\widetilde B_P=a\widetilde B_P$ and hence $b\widetilde{B}_P\cap{B_{\p}}=aB_{\p}$. Therefore $D_{\p}(B_{\p})B_{\p}\subseteq aB_{\p}$, contradicting the fact that $D_{\p}$ is irreducible. So, in this case we have $\widetilde D_P$ is irreducible.  

\smallskip\noindent
{\bf Case 2:} Suppose, $rf\in\widetilde B\setminus C$. If $g\in C^*~(=R^*)$, then considering $r,b,f,a,g$ as a polynomials in $X_1,\dots,X_{n_1}, Y_1,\dots,Y_{n_2}$ with coefficients  from  $B$ and comparing the total degree we can see that $b\in B$, which gives $b\widetilde B_P\cap B_\p=bB_\p$. Then, as before we arrive at a state of contradiction, and hence $\widetilde D_P$ is irreducible. So, without loss of generality we can assume that $g$ is a non-unit of $C$. Then, there exists $Q\in\Spec{C}$ such that $g\in Q$. Since $(B\otimes_RC)\otimes_Ck(Q)=k(Q)^{[2]}$, the equation $brf=ag$ implies either $b\in QC_Q$ or $rf\in QC_Q$. Since $R$ is inert in $B$, we have $C$ is inert in $\widetilde B$ and hence $QC_Q\cap \widetilde B=Q$. So, $b\in Q$ (as  $rf\in\widetilde B\setminus C$). Since $b~(\in C)$ is a non-unit of $\widetilde B$, it follows that $b\in P$ and hence $ag\in P$, contradicting the fact $a\in R\setminus\p$ and $g\in C\setminus P$.

\smallskip\noindent Thus, in any situation, $\widetilde D_P$ is irreducible. Therefore, by Proposition \ref{Prop_Triv-A2-Fib_Sym-Ker} we have $\grade{\widetilde{D}(\widetilde B)\widetilde B,\widetilde B}$ is either $2$ or $\infty$. Now, (I) follows from Lemma \ref{Lem_equality of grade under ff}. 
\end{proof}

\begin{rem}\label{rem_triviality A2-fib}
In Theorem \ref{Thm_A2-Fib_A1-Ker}, when $\grade{D(B),B}$ is $\infty$, it has been proved in \citep[Theorem 4.4]{Babu-Das_Struct_A2-fib_FPF-LND} that $B=A^{[1]}$. Moreover, if $R$ is seminormal, then $B=({\rm Sym}_{R}(N))^{[1]}$ for some rank one projective $R$-module $N$.
\end{rem}

\indent In \citep[Theorem 3.8]{Dutta-Gupta-Lahiri_separable-A2-A3-forms}, it has been proved  that for an $R$-lnd $D$ on an $\A{2}$-form $B$ over a ring $R$ containing $\bQ$, if $\grade{D(B)B,B}=\infty$, then $\Ker{D}= {\rm Sym}_{R}(N)$ for some rank one projective $R$-module $N$ and $B=A^{[1]}$. The next corollary explores the case when $\grade{D(B)B,B}=2$.
\begin{cor}\label{Cor_A2-forms_Sym-Ker}
Let $k$ be a field of characteristic zero, $R$ a Noetherian $k$-domain  and $B$ an $\A{2}$-form over $R$. Suppose, $D\in{\rm LND}_{R}(B)$ and $A={\rm Ker}(D)$. Then the following are equivalent.
		
\begin{enumerate}
\item [\rm (I)] ${\rm grade}\big{(}D(B)B,B\big{)}$ is either $2$ or $\infty$.
			
\item [\rm (II)] $A={\rm Sym}_{R}{(N)}$ for some rank one projective $R$-module $N$ and  $D_{\p}$ is irreducible for each $\p\in {\rm Spec}(R)$. 
\end{enumerate}

\end{cor}
\begin{proof}
 Since $\A{2}$-forms over a ring containing field of characteristic zero are $\A{2}$-fibrations over the corresponding ring (See \cite[Lemma 3.6]{Dutta-Gupta-Lahiri_separable-A2-A3-forms}), (II)$\implies$(I) follows from Theorem \ref{Thm_A2-Fib_A1-Ker}.
 
 \medskip
 \noindent
 (I)$\implies$(II): Suppose, (I) holds. By Theorem \ref{Thm_A2-Fib_A1-Ker}, $D_\p$ is irreducible for each $\p\in\Spec{R}$. Let $L$ be a finite extension of $k$ such that $B\otimes_kL=(R\otimes_kL)^{[2]}$. Set $\widetilde{R}=R\otimes_kL, \widetilde{A}=A\otimes_kL, \widetilde{B}=B\otimes_kL$ and $\widetilde{D}=D\otimes_k1_L$. Then, $\widetilde{B}=\widetilde{R}^{[2]}$ and $\grade{D(\widetilde{B})\widetilde{B},\widetilde{B}}$ is $2$ or $\infty$. By Theorem \ref{BD_LND-Characterisation-Result_Grade}, we have $\widetilde{A}=\widetilde{R}^{[1]}$ and therefore by Theorem \ref{AKD1}, $A=\Sym{R}{N}$ for some rank one projective $R$-module $N$.
\end{proof}

In \cite{Babu-Das-Lokhande_Rank-and-Rigidity}, Babu-Das-Lokhande introduced the concept of residual rank and residual-variable rank of lnds on affine fibrations.
When $\grade{D(B)B,B}$ is $2$ or $\infty$, the next result provide some sufficient conditions, in terms of residual rank (or residual variable rank) of an $R$-lnd $D$ on an $\A{n}$-fibration over $R$, for $\Ker{D}$ to be $\A{n-1}$-fibration over $R$.

\begin{cor}\label{Cor_resrank}
	Let $R$ be a Noetherian domain containing $\bQ$, $B$ an $\A{n}$-fibration over $R$ and $D \in {\rm LND}_R(B)$ with ${\rm grade}(D(B)B, B)\in\{2, \ \infty\}$. Then the following statements hold.
	\begin{enumerate}
		\item[\rm(1)]If ${\rm Res-Rk}(D) =2$, then ${\rm Ker}(D)$ is an $\A{n-1}$-fibration over $R$.
		\item[\rm(2)]If ${\rm ResVar-Rk}(D) =2$, then ${\rm Ker}(D)$ is an $\A{1}$-fibration over $R^{[n-2]}$.
	\end{enumerate}
	
\end{cor} 

\begin{proof}
	{(1):} Assume that $\ResRk{D} =2$. Then, there exists an $(n, 2)$-residual system $(R, C, B)$ such that $C \subseteq \Ker{D}$. Therefore, $B$ is an $\A{2}$-fibration over $C$ and $D \in \LND{C}{B}$. By Theorem \ref{Thm_A2-Fib_A1-Ker}, we have $\Ker{D}$ is an $\A{1}$-fibration over $C$ and hence for each $\p \in \Spec{R}$, $\Ker{D} \otimes_{R} k(\p)$ is an $\A{1}$-fibration over $C \otimes_{R} k(\p)$. Since $C$ is an $\A{n-2}$-fibration over $R$, we have $C \otimes_{R} k(\p)=k(\p)^{[n-2]}$, a UFD. Since rank one projective modules over UFDs are free, by Corollary \ref{Cor_A1-fib_Triv_DiffMod_Extended}, $\Ker{D} \otimes_{R} k(\p) = \big(C \otimes_{R} k(\p)\big)^{[1]} = k(\p)^{[n-1]}$. Faithful flatness and finite generation of $\Ker{D}$ over $R$ is clear.
	
\medskip
\noindent
{(2):} Assume that $\ResVarRk{D}=2$. Then, there exists an $(n, 2)$-residual system $(R, C, B)$ such that $C \subseteq \Ker{D}$ and $C=R^{[n-2]}$. Repeating the arguments in (1), it is easy to see that $\Ker{D}$ is an $\A{1}$-fibration over $C = R^{[n-2]}$.
\end{proof}

\section{Structure of the kernel over Noetherian normal domain}\label{Sec_Normal}

\noindent
In this section  we answer part (b) of Question \ref{Qtn_main}. First, we prove the following proposition which generalizes Proposition 3.3 of \citep{Bhatwadekar-Dutta_LND}. Our proofs are highly inspired by the proof of Proposition 3.3 in \citep{Bhatwadekar-Dutta_LND}.
\begin{prop} \label{Prop_Inert-subring-of-affine-fib}
Let $R$ be a Noetherian normal domain, $B$ an $\A{r}$-fibration over $R~(r\geqslant 1)$ and $A$ an inert $R$-subalgebra of $B$ with $\trdeg{B}{R}=1$. Then, $A$ has the structure of a graded ring $\bigoplus_{i\ge 0} A_i$ with $A_0 = R$ and for each $i \geqslant 1$, $A_i$ is a finite reflexive $R$-module of rank one. In fact, when $R$ is not a field, then there exists an unmixed height one ideal $I$ of $R$ such that $A$ is isomorphic to the symbolic Rees algebra $\bigoplus_{n\ge 0} I^{(n)}T^{n}$.
\end{prop}

\begin{proof} By Theorem \ref{Asanuma_struct-fib-th}, there exists $m \in \mathbb{N}$ such that $B$ is a retract of $C = R^{[m]}$.

\medskip\noindent
{\bf Case 1.} If $R$ is a field, then $B=R^{[r]}$ and hence by a well-known result of L\"uroth we have $A=R^{[1]}$. 

\medskip\noindent
{\bf Case 2.} Assume that ${\rm dim}(R)=1$ and $\m$ is an arbitrary maximal ideal of $R$. Then, $R_\m$ is a discrete valuation ring and hence $C_\m$ is a UFD. Since $B_\m$ is a retract of $C_\m$ and $A_\m$ is an inert subring of $B_\m$, we see that $A_\m$ is a UFD. Since ${\rm tr.deg}_{R_\m}(A_\m)=1$, by Theorem \ref{AEH_Luroth} we have $A_\m= (R_\m)^{[1]}$. Now, by Proposition \ref{Giral} and Theorem \ref{Onoda_FG}(II), $A$ is finitely generated over $R$. By Theorem \ref{BCW_Sym}, $A=\Sym{R}{M}$ for some projective $R$-module $M$ of rank one. Since $M$ is projective of rank one, $M$ is isomorphic to an invertible ideal of $R$.
 
\medskip\noindent
{\bf Case 3.} Assume that ${\rm dim}(R)\geqslant 2$.
 Since $A\subseteq C~(=R^{[m]})$, we have $JA \cap R = J$ for every ideal $J$ of $R$. Therefore, by Lemma \ref{Lem_BD-A1-patch}, it is enough to show that there exists an $R$-sequence $\{x,\ y\}$ such that $A_x = R_x^{[1]}$, $A_y = R_y^{[1]}$ and $A = A_x \cap A_y$.	
Let $T = R \setminus \{0\}$ and $K=T^{-1}R$. Since $T^{-1}A$ is an inert subring of $T^{-1}B ~(= B \otimes_R K= K^{[2]})$ of transcendence degree one over $K$, by Theorem \ref{AEH_Luroth} we have $T^{-1}A =K^{[1]}$. Since $A\subseteq C~(=R^{[m]})$, by Proposition \ref{Giral} and Theorem \ref{Onoda_FG}(I), there exists $t\in T$ such that $A_t$ is finitely generated as an $R$-algebra and hence using the equality $T^{-1}A =K^{[1]}$, we can choose $x\in T$ such that $A_x = R_x^{[1]}$. If $x \in R^*$, then $A = R^{[1]}$, and we are done. So, we assume that $x \notin R^*$. Let  $P_1, \cdots, P_{\ell}$ be the associated primes of $R/xR$. Since $R$ is a Noetherian normal domain, $\Ht{P_i} = 1$ for all $i = 1, 2, \cdots, \ell$. Set $ S:=R\backslash\left( \bigcup_{i =1}^{\ell} P_i \right)$. Then, $S^{-1}R$ is a PID and hence $ S^{-1}C$ is a UFD. Since $S^{-1}B$ is a retract of $S^{-1}C$ and $S^{-1}A$ is an inert subring of $S^{-1}B$, we see that $S^{-1}A$ is a UFD. Since $\trdeg{S^{-1}A}{S^{-1}R}=1$, by Theorem \ref{AEH_Luroth} we have $ S^{-1}A= (S^{-1}R)^{[1]}$ and hence applying Proposition \ref{Giral} and Theorem \ref{Onoda_FG}(I) once again we choose $y\in S$ such that $A_y=R_y^{[1]}$.
By construction, the pair $x,y$ form an $R$-regular sequence. We now show that $A = A_x \cap A_y$. Let $c=a/x^j =b/y^l \in A_x \cap A_y$; $a,b \in A$; $j,l\geqslant 1$. Since $x, y$ form an $R$-regular sequence and $B$ is faithfully flat over $R$, we see that $x,y$ form a $B$-regular sequence. Therefore, from the equation $y^la = x^jb$ it follows that $a \in  x^jB$. Since $A$ is an inert subring of $B$, we have $a\in x^jA~(=x^jB\cap A)$. Therefore, $c = a/x^j \in A$, which shows that $A = A_x \cap A_y$. This completes the proof.
\end{proof}

\begin{thm}\label{Thm_kernel_Noeth_normal}
Let $R$ be a Noetherian normal domain containing $\bQ$, $B$ an $\A{2}$-fibration over $R$ and $D\in\LND{R}{B}\setminus\{0\}$. Then, $\Ker{D}$ has the structure of a graded ring $\bigoplus_{i\ge 0} A_i$ with $A_0 = R$ and for each $i \geqslant 1$, $A_i$ is a finite reflexive $R$-module of rank one. In fact, when $R$ is not a field, then there exists an ideal $I$ of unmixed height one in $R$ such that $A$ is isomorphic to the symbolic Rees algebra $\bigoplus_{n\ge 0} I^{(n)}T^{n}$.

\smallskip\noindent
Conversely, let $R$ be as above and let $I$ be an unmixed ideal of height one in $R$. Let $A$ be the symbolic Rees algebra
$\bigoplus_{n\geqslant 0}I^{(n)}T^n$. Then, there there exists an $\A{2}$-fibration $C$ over $R$ and lnd $D \in \LND{R}{C}$ such that $\Ker{D}$ is isomorphic to $A$ as a graded $R$-algebra. In particular $A$ can be embedded as an inert subring of $C$.
\end{thm}
\begin{proof}
Since $\Ker{D}$ is an inert subring of $B$ with ${\rm tr.deg}_{R}(B)=1$, the first part of the theorem follows from Proposition \ref{Prop_Inert-subring-of-affine-fib}.

\smallskip\noindent
For the converse part, let $A$ be the symbolic Rees algebra $\bigoplus_{n\geqslant 0}I^{(n)}T^n$, where $I$ is an unmixed ideal of height one in $R$. By \cite[Theorem 3.5]{Bhatwadekar-Dutta_LND}, there exists $D \in \LND{R}{R[X,Y]}$ such that $\Ker{D}$ is isomorphic to $A$ as a graded $R$-algebra, and therefore, $A$ can be embedded as an inert subring of $R[X,Y]$. Now one can see that the proof is done if we set $C: =R[X,Y]$.
\end{proof}

\begin{rem} ~
\begin{enumerate}

\item Since Proposition \ref{Prop_Inert-subring-of-affine-fib} does not depend upon the characteristic of $R$, Theorem \ref{Thm_kernel_Noeth_normal} can be proved for exponential maps on $B$.

\item If $R$ is a Noetherian normal domain, $B$ is an $\A{2}$-fibration over $R$ and $A$ is an inert subring of $B$ with transcendence degree one over $R$, then we have shown that the structure of $A$ is exactly same as the case of $B=R^{[2]}$. This raises the question that whether any $\A{2}$-fibration over $R$ is locally polynomial algebra over $R$? If $R$ contains $\bQ$, then this question is OPEN till now. However, when $R$ is a normal domain not containing $\bQ$, there exist $\A{2}$-fibrations over $R$ which are not locally polynomial (see \cite{Asanuma_fibre_ring}).
  
\end{enumerate}
\end{rem}
\section{Examples}\label{Sec_Examples}
\noindent
In this section our main aim is to discuss examples (see Example \ref{Ex_red1} and Example \ref{Ex_red2}) to demonstrate that the equivalent condition (II) in Theorem \ref{Thm_A2-Fib_A1-Ker} can't be reduced further. However, we start with the following example where we have constructed two non-trivial $\A{2}$-fibrations $B$ and $C$ over $R=k[X^2,X^3]$ having  kernels isomorphic to a non-trivial $\A{1}$-fibration over $R$ with grade $2$ and $\infty$ respectively.

\begin{ex}\label{Ex_grade}
 Let $k$ be a field of characteristic zero, $R=k[X^2, X^3]$, $A=R[Y+ XY^2] + X^2k[X,Y]$. In \cite[Example 1, Section 4]{Yanik}, it has been proved that $A$ is a retraction of a polynomial algebra over $R$ and hence is flat as an $R$-module and finitely generated as an $R$-algebra. 
 Suppose, $\p_0=(X^2,X^3)R$. Then $A\otimes_R\dfrac{R}{\p_0}=\dfrac{R}{\p_0}[\overline{Y+XY^2}]$. Suppose, $\p$ is any prime ideal of $R$ other than $\p_0$. Since $X\in R_\p$, one can see that $A\otimes_RR_\p=R_\p[Y]$. Therefore, $A$ is a non-trivial $\A{1}$-fibration over $R$.
 
\medskip\noindent
Let $B=A[Z]$. Clearly, $B$ is a non-trivial $\A{2}$-fibration over $R$. Let $D=(\frac{d}{dZ})_{A}\in\LND{R}{B}$. Then $\Ker{D}=A$ and $\grade{D(B)B,B}$ is $\infty$.
 
\medskip\noindent 
Let $C=A[Z+XZ^2]+X^2k[X,Y,Z]$. Then $C~(=B\otimes_RA)$ is a non-trivial $\A{1}$-fibration over $A$ and hence a non-trivial $\A{2}$-fibration over $R$. Let $\widetilde{D'}=X^2\dfrac{\partial}{\partial Z}\in \LND{R}{k[X,Y,Z]}$. Then $\widetilde{D'}(C)\subseteq C$. Define $D'=\widetilde{D'}|_{C}$. Then $D'\in\LND{R}{C}$. Clearly, $A\subseteq\Ker{D'}$. Since $C$ is an $\A{1}$-fibration over $A$, it is easy to check that $A$ is inert in $C$ and hence $\Ker{D'}=A$. Clearly, $D'(C)C=(X^4,X^2+2X^3Z)C$. Since $X\notin B$, we see that $X^4,X^2+2X^3Z$ form a $C$-regular sequnce; and hence by Lemma \ref{Lem_2genIdeal_seq-grade}, $\grade{{D'}(C)C,C}$ is $2$. 
\end{ex}

Now, we give an example to show that in the equivalent condition (II) of  Theorem \ref{Thm_A2-Fib_A1-Ker}, the condition ``$D_\p$ is irreducible for each prime ideal $\p$ of $R$'' is not redundant. 
\begin{ex}\label{Ex_red1}
 Let $R,A,C$ as in Example \ref{Ex_grade} and $\widetilde{D}=X^4\dfrac{\partial}{\partial Z}\in \LND{R}{k[X,Y,Z]}$. Then $\widetilde{D}(C)\subseteq C$. If we define $D=\widetilde{D}|_{C}$, then $D\in\LND{R}{C}$ and $\Ker{D}=A$. It is easy to see that $D(C)C=(X^6,X^4+2X^5Z)C\subseteq X^2C$. Since $X^4(X^4+2X^5Z)\in X^6C$, $\grade{{D}(C)C}=1$. Here, we note that if we take  $\p_0=(X^2,X^3)\in\Spec{R}$, then $D_{\p_0}(C_{\p_0})\subseteq X^2C_{\p_0}$ and hence $D_{\p_0}$ is reducible. We also note that for any other prime ideal $\p$ of $R$, $X$ being a unit in $R_\p$, $D_\p$ is irreducible.
\end{ex}

Next, we quote an example of Bhatwadekar and Dutta (see \cite [Example 3.11]{Bhatwadekar-Dutta_LND}) to show that in the equivalent condition (II) of Theorem \ref{Thm_A2-Fib_A1-Ker}, the condition ``$A$ is an $\A{1}$-fibration'' is not redundant.
\begin{ex}\label{Ex_red2}
 Let $R=\bR+(X)\bC[[X]]$, $S=\bC[[X]]$ and $B=R[Y,Z]$. Define a $S$-lnd $\widetilde{D}$ on $S[Y,Z]$ by setting $\widetilde{D}(Y)=iX$ and $\widetilde{D}(Z)=-X$. Then $D:=\widetilde{D}|_B\in\LND{R}{B}$. It has been shown in \cite [Example 3.11]{Bhatwadekar-Dutta_LND}) that $A=\Ker{D}$ is not a finitely generated as an $R$-algebra.
Clearly, $\grade{D(B)B,B}=1$ and
$D_\p$ is irreducible for each prime ideal of $R$. In fact, if $\p=(0)$, then $D_\p$ is an $R$-lnd with a slice and if $\p=(X)$, then $D_\p$ is irreducible as $i\notin R$.
\end{ex}

\section*{Acknowledgment} 
\noindent
The authors thank Neena Gupta for her valuable comments. The second author acknowledges SERB, Govt. of India for their MATRICS grant under the file number MTR/2022/000247.

	\bibliographystyle{alpha}
	\normalem
	\bibliography{reference}
\end{document}